\documentclass{article}
\usepackage[a4paper]{geometry}
\usepackage{amssymb}
\usepackage{amsmath}
\usepackage{rotating}
\usepackage{pdflscape}
\usepackage{longtable}
\usepackage[table]{xcolor}
\usepackage{epsfig}
\usepackage {graphicx}
\usepackage{caption}

\definecolor{light-gray}{gray}{0.95}
\numberwithin{equation}{section}
\newcommand{\beq}{\begin{equation}}
\newcommand{\eeq}{\end{equation}}
\newtheorem{Thm}{Theorem}[section]

\newtheorem{lemma}[Thm]{Lemma}

\newtheorem{theorem}[Thm]{Theorem}

\newcommand{\ds}[0]{\,\dot{+}\,}

\newcommand{\R}[0]{\mathbb{R}}

\newcommand{\N}[0]{\mathbb{N}}

\newcommand{\s}[0]{\subset}
\newcommand{\mc}[1]{\mathcal{#1}}
\newcommand{\x}[0]{\times}

\begin{document}
\title{Automorphisms of real Lie algebras of dimension five or less}
\author{David J Fisher\thanks{D.J.Fisher@surrey.ac.uk}, Robert J Gray\thanks{R.J.Gray@surrey.ac.uk}, Peter E Hydon\thanks{P.Hydon@surrey.ac.uk}\\ \\Department of Mathematics, University of Surrey, Guildford GU2 7XH, UK}
\maketitle

\begin{abstract}
The Lie algebra version of the Krull-Schmidt Theorem is formulated and proved. This leads to a method for constructing the automorphisms of a direct sum of Lie algebras from the automorphisms of its indecomposable components.  For finite-dimensional Lie algebras, there is a well-known algorithm for finding such components, so the theorem considerably simplifies the problem of classifying the automorphism groups. We illustrate this by classifying the automorphisms of all indecomposable real Lie algebras of dimension five or less. Our results are presented very concisely, in tabular form.
\end{abstract}

\section{Introduction}
Automorphisms of Lie algebras are of fundamental importance in various branches of mathematics and physics.
For example, nilpotent Lie algebras of a given dimension can be classified by treating each one as a central extension of a lower-dimensional nilpotent Lie algebra, $\mc{L}$. The automorphisms of $\mc{L}$ are used to sort the set of extensions into equivalence classes \cite{CdeGS,SkS}, greatly simplifying the classification problem. Similarly, if one knows the group of automorphisms of each Lie algebra of dimension $n$ over a field of characteristic zero, one can (in principle) classify all solvable Lie algebras of dimension $n+1$ over the same field. However, the complexity of the computations increases very rapidly with $n$. An alternative approach, due to Patera \& Zassenhaus \cite{PaZ}, is to construct each solvable Lie algebra of a given dimension from an associated nilpotent Lie algebra of the same dimension; once again, automorphisms are used to simplify the calculations.

Lie algebra automorphisms are useful in cosmology; Terzis \& Christodoulakis \cite{PaC} used them to find systematically a complete set of vacuum solutions to Einstein's field equations for Bianchi-type spacetime geometries. Similarly, the problem of classifying four-dimensional homogeneous spaces for $4+1$-dimensional spacetimes led Christodoulakis \textit{et al.} to classify the automorphisms of every four-dimensional real Lie algebra \cite{CPD}.
 
Automorphisms are also essential for constructing the discrete symmetries of a given differential equation. Such symmetries are used in equivariant bifurcation theory, the construction of invariant solutions and the simplification of numerical schemes. Given a differential equation, the symmetry condition is a highly-coupled nonlinear system; one cannot usually determine discrete symmetries by attempting to solve the symmetry condition directly.  By contrast, it is often quite easy to solve the linearized symmetry condition, which yields generators of Lie point, contact or generalized symmetries. Once a nontrivial Lie algebra $\mathcal{L}$ of symmetry generators is known, the problem of finding all discrete symmetries can be tackled indirectly, as follows \cite{Hy1,Hy2}. The adjoint action of each symmetry induces an automorphism of the Lie algebra, so the first step is to determine all automorphisms of $\mc{L}$. To simplify the problem, it is helpful to remove some parameters by applying suitable inner automorphisms. The remaining automorphisms produce a set of constraints that determine the discrete symmetries. Typically, it is easy to solve the symmetry condition subject to the constraints. Indeed, the hardest part of the calculation is the classification of all automorphisms of $\mathcal{L}$. For this reason, it
is valuable to have look-up tables that list the automorphisms of various types of real Lie algebras.

At present, the largest such
tables are those of Christodoulakis \textit{et al.} \cite{CPD} and Laine-Pearson \& Hydon \cite{LPH}.
The classification by Christodoulakis \textit{et al.} deals with all four-dimensional real Lie algebras; it is based on
the comprehensive list of real Lie algebras in Patera \textit{et al.} \cite{PSWH}, which is adapted from the work of Mubarakzyanov \cite{Mub1,Mub2,Mub3}. Earlier, Harvey
classified the automorphisms of three-dimensional Lie algebras \cite{Har}; the increase from three to four dimensions substantially increases the complexity
of the classification.
Laine-Pearson \& Hydon's table is based on the
classification by Gonz\'{a}lez-L\'{o}pez \textit{et al.} of all real Lie algebras that occur for point symmetries of scalar ordinary differential equations \cite{GLKO}.
This classification includes some families of Lie algebras of arbitrary dimension, but it excludes many Lie algebras that occur for systems of
ordinary differential equations and for partial differential equations.

Both of the above classifications are lengthy, so one might expect that a classification of automorphisms of five-dimensional Lie algebras would be too long for publication. However, the results can be greatly simplified by removing various redundancies.
A Lie algebra $\mathcal{L}$ is \textit{decomposable} if it is a direct sum of proper ideals; if these ideals are indecomposable, the decomposition is \textit{irreducible}.
Given any decomposable finite-dimensional Lie algebra, one can construct an irreducible decomposition by using an algorithm due to Rand \textit{et al.} \cite{RWZ} (see also \cite{deGraaf}).

The expression of a group as an internal direct product of a finite number of indecomposable normal subgroups is known to be unique, up to isomorphism and the order of the factors, provided that the ascending and descending chain conditions hold on normal subgroups. This is the Krull-Schmidt Theorem (see \cite{Hungerford}), which also establishes an exchange property for such groups:
if
\[
G = H_1 \x H_2 \x \cdots \x H_r = J_1 \x J_2 \x \cdots \x J_s
\]
then $r = s$ and the factors can be labelled such that, for each $i, \, H_i \cong J_i$ and
\[
G = H_1 \x \cdots \x H_i \x J_{i+1} \x \cdots \x J_r. 
\]
A corresponding theorem holds for a module over a ring, expressed as a direct sum of indecomposable submodules, provided that the chain conditions hold for submodules (\cite{Jacobson1}).

We have not found a similar result for Lie algebras in the literature.  Theorem 2.8 of \cite{RWZ} summarizes the decomposition algorithm, but Theorem 2.8 (e) seems to imply (incorrectly) that the decomposition is unique when there is no central component.  We show that uniqueness holds if the algebra is centreless or equals its own derived algebra.  Otherwise, multiple decompositions are possible, which we characterize.  Our proof is modelled on those for groups and modules. 

Having proved the Krull-Schmidt Theorem for Lie algebras, we introduce a method for constructing the automorphisms of a decomposable Lie algebra from the automorphisms of its indecomposable components.  Armed with such a method, it is then sufficient to classify the automorphisms of indecomposable Lie algebras.

Finally, we classify all automorphisms of five-dimensional indecomposable real Lie algebras and present the results in a compact look-up table. This allows the user to identify inner and outer automorphisms, distinguishing between exponentiated derivations and discrete automorphisms.
For completeness, we also present the automorphisms of lower-dimensional indecomposable real Lie algebras in the same form.

\textbf{Note}: As our main motivation for producing look-up tables of Lie algebra automorphisms is their use for finding discrete symmetries, we use the standard convention for symmetry analysis, which is that the Lie algebra of a symmetry Lie group $G$ is the vector space of \textit{right}-invariant vector fields on $G$ \cite{Olver}. For applications where the Lie algebra is regarded as the set of left-invariant vector fields on $G$, the automorphism matrices that arise from our tables must be transposed.

\section{Determining equations for automorphisms}
Given a basis $\{X_1,\cdots,X_R\}$ of a real finite-dimensional Lie algebra $\mathcal{L}$, the structure constants $c_{ij}^k$ are determined by the commutators of the basis elements,
\[
 \big[X_i, X_j\big]=c_{ij}^kX_k,\qquad 1\le i<j\le R.
\]
(Here and throughout, the Einstein summation convention applies.)
An automorphism of $\mathcal{L}$ is a bijective linear map $\Phi:\mathcal{L}\mapsto \mathcal{L}$ such that for all $X,Y\in \mathcal{L}$,
\beq\label{aut1}
 \Phi\big([X,Y]\big)=\big[\Phi(X), \Phi(Y)\big].
\eeq
In terms of the given basis, the automorphism $\Phi$  is represented by a
real-valued matrix $B=\big(b_i^l\big)$, where the lower (upper) index is the row (column) number:
\[
 \Phi(X_i) = b_i^lX_l,\qquad i=1,\dots,R.
\]
Then the automorphism condition (\ref{aut1}) amounts to
\beq\label{AC}
c_{lm}^nb_i^lb_j^m = c_{ij}^kb_k^n; \qquad \text{det}(B)\neq 0.
\eeq
The equations with $i \geq j$ are equivalent to those with $i<j$, so (\ref{AC}) is a system of up to $R^2(R-1)/2+1$ determining equations for the automorphism group $\mathrm{Aut}(\mathcal{L})$. Typically, the complexity of this system grows rapidly with $R$, though it is determined mainly by the number of structure constants that are nonzero.

It is usually possible to spot some zero sub-blocks of $B$ without having to solve any equations.
For example, any automorphism will map the elements of the derived series of $\mc{L}$ to itself.
The same is true for the upper and lower central series of $\mc{L}$.
Furthermore, \eqref{AC} yields the following simple necessary conditions for any automorphism (see \cite{CPD}).
\[
c_{ln}^nb_j^l=c_{ln}^n\delta^l_j, \qquad c_{lk}^nc_{mn}^kb_i^lb_j^m=c_{ik}^nc_{jn}^k;
\]
here $\delta^l_j$ is the Kronecker delta.
The first of these identities shows that if $\mathrm{tr\{ad}(X_l)\}$ (which is $-c_{ln}^n$) is nonzero then $1$ is an eigenvalue of every automorphism matrix; the second expresses the fact that every automorphism preserves the Killing form.

As $\mathcal{L}$ is finite-dimensional, the group of automorphisms, $\mathrm{Aut}(\mc{L})$, is a Lie group. Our aim is to obtain the matrix representation of this group by solving the determining equations (\ref{AC}). Typically, some of these equations are very easy to solve, but there will often be a remaining set of highly-coupled conditions.

It is helpful to simplify these by using appropriate inner automorphisms, $\mathrm{Ad}_X$. The adjoint action of each basis vector, $X_j$, generates a one-parameter Lie group of inner automorphisms whose matrix representation is
\beq\label{adxj}
A_j(\varepsilon_j)=\exp\{\varepsilon C(j)\},\qquad \text{where}\quad \big(C(j)\big)^k_i=c_{ij}^k.
\eeq
(Here $\exp$ is the matrix exponential.) In the matrix representation, every inner automorphism is a finite product of the matrices \eqref{adxj}. Indeed, every inner automorphism that is sufficiently close to the identity can be expressed as a product in which each $A_j$ is used \textit{once}. We can use any convenient ordering of these matrices, but the values of the parameters $\varepsilon_j$ will depend on the ordering. (Generally speaking, this result is not true globally;  Weyl reflections are inner automorphisms that cannot be written in this way.) Nevertheless, the local result is useful, as follows.

Suppose that we have solved as many of the determining equations \eqref{AC} as possible and therefore $B$ is partly determined. Let
\beq\label{tb}
\tilde{B}=A_1(\tilde{\varepsilon}_1)A_2(\tilde{\varepsilon}_2)\cdots A_R(\tilde{\varepsilon}_R)B.
\eeq
Clearly, $B$ satisfies the remaining determining equations if and only if $\tilde{B}$ does. Thus we choose the parameters $\tilde{\varepsilon}_j$ so as to create as many zero entries in $\tilde{B}$ as possible, and we use any $A_j$ that are diagonal to scale one or more nonzero entries to $\pm 1$. It is not necessary to use the ordering in \eqref{tb}; any convenient ordering will do, so the simplification can proceed one step at a time.

This approach is highly effective. Once the remaining determining equations have been solved and we know $\tilde{B}$, we can reconstruct the whole automorphism group by writing
\beq\label{bt}
B=A_1(\varepsilon_1)A_2(\varepsilon_2)\cdots A_R(\varepsilon_R)\tilde{B},
\eeq
where now each $\varepsilon_j$ is a free variable. 

One further major simplification is possible for any finite-dimensional Lie algebra $\mc{L}$ that is a direct sum of proper ideals, $\mc{M}_i$: one can construct $\mathrm{Aut}(\mathcal{L})$ from each $\mathrm{Aut}(\mathcal{M}_i)$. This result uses the Krull-Schmidt Theorem for Lie algebras, which we state and prove in the next section. 

\section{Automorphisms of direct sums of Lie algebras}

Let $\mc{L}$ be a Lie algebra that is a direct sum of proper ideals:
\beq\label{dsd}
\mc{L}=\mc{M}_1 \oplus \cdots \oplus \mc{M}_r.
\eeq
In other words, $\mc{L}$ is a vector space direct sum (denoted $\ds$),
\[
\mc{L}=\mc{M}_1 \ds \cdots \ds \mc{M}_r,\qquad \mc{M}_i\cap\mc{M}_j= \{0\},
\]
and $[\mc{M}_i, \mc{M}_j] = \{0\}$ for $i\neq j$, because each $\mc{M}_i$ is an ideal. We call a direct sum decomposition of $\mc{L}$ \textit{irreducible} if each $\mc{M}_i$ is {\em indecomposable}, that is, if no $\mc{M}_i$ can be expressed as a direct sum of proper ideals.  We denote the centre of $\mc{L}$ by $\mc{Z}(\mc{L})$ and the derived algebra, $[\mc{L},\mc{L}]$, by $\mc{L}'$.  A direct summand
$\mc{M}_i$ is a {\em central component} of $\mc{L}$ if $\mc{M}_i \s \mc{Z}(\mc{L})$. Clearly, each central component is abelian, so if it is indecomposable, it has dimension 1.

Given any direct sum decomposition (\ref{dsd}) of $\mc{L}$, each $X \in \mc{L}$ has a unique representation:
\[
X = X_1 + \cdots + X_r,\qquad X_i\in \mc{M}_i.
\]
The projection $\pi_i$ is the endomorphism of $\mc{L}$ defined by $\pi_i(X) = X_i$, so $\mc{M}_i = \pi_i(\mc{L})$.
Clearly,
\[
{\pi_i}^2 = \pi_i, \qquad \pi_i\pi_j = 0\text{ if }i\ne j,\qquad\pi_1 + \cdots + \pi_r = \iota,
\]
where $\iota$ is the identity map on $\mc{L}$. Let $\overline{\mc{M}}_i = (\iota - \pi_i)(\mc{L})$; this is a direct sum of all $\mc{M}_j$ for $j \ne i$.

An endomorphism $\phi: \mc{L} \to \mc{L}$ is {\em normal} if $\phi$ commutes with each $\mathrm{ad}_X$: 
\[
\phi([X,Y]) = [\phi(X), Y] = [X, \phi(Y)],\qquad X, Y \in \mc{L}.
\] 
In particular, the projection maps $\pi_i$ and $\iota-\pi_i$ are normal endomorphisms of $\mc{L}$.  Any composite of normal endomorphisms is normal. 
If $\phi$ is an endomorphism of $\mc{L}$ and
\beq\label{fitcon}
\phi^2(\mc{L}) = \phi(\mc{L}),\qquad \text{Ker}(\phi^2) = \text{Ker}(\phi),
\eeq
then, by Fitting's Lemma (\cite{Jacobson2}), $\mc{L} = \text{Ker}(\phi) \ds \phi(\mc{L})$.  If $\phi$ is normal then $\phi(\mc{L})$ is an ideal of $\mc{L}$, so when the conditions (\ref{fitcon}) hold, there is a direct sum decomposition, $\mc{L} = \text{Ker}(\phi) \oplus \phi(\mc{L})$.

A Lie algebra $\mc{L}$ is {\em Artinian} or {\em Noetherian} if it satisfies the descending or ascending chain condition respectively on ideals.  The ascending chain condition guarantees that any direct-sum decomposition has finitely many components. If $\mc{L}$ is finite-dimensional, it is both Artinian and Noetherian. However, the results that we prove below apply equally to infinite-dimensional Lie algebras that satisfy both chain conditions.

\begin{lemma}\label{L1}
Let $\mc{M}$ be an indecomposable Lie algebra which is both Artinian and Noetherian. Then every normal endomorphism of $\mc{M}$ is either bijective or nilpotent.
\end{lemma}

\noindent{\bf Proof}\, Let $\phi : \mc{M} \to \mc{M}$ be a normal endomorphism, so $\phi(\mc{M})$ is an ideal of $\mc{M}$.  As both chain conditions hold for ideals, there exists $k \in \N$ such that for all $\ell \ge k, \, \phi^k(\mc{M}) = \phi^{\ell}(\mc{M})$ and $\text{Ker}(\phi^k) = \text{Ker}(\phi^{\ell})$.  Taking $\ell = 2k$, Fitting's Lemma gives $\mc{M} = \text{Ker}(\phi^k) \oplus \phi^k(\mc{M})$.

As $\mc{M}$ is indecomposable, either $\mc{M} = \text{Ker}(\phi^k)$, so $\phi$ is nilpotent, or $\mc{M} = \phi^k(\mc{M})$ and $\text{Ker}(\phi^k) = \{0\}$, in which case $\phi^k$ is bijective and hence so is $\phi$. \hfill $\Box$ \\

For the rest of this section, we shall assume that $\mc{L}$ is both Artinian and Noetherian. Suppose that $\mc{M}_1 \oplus \cdots \oplus \mc{M}_r$ and $\mc{N}_1 \oplus \cdots \oplus \mc{N}_s$ are irreducible decompositions of $\mc{L}$, with projections $\pi_i$ (onto $\mc{M}_i$) and $\psi_j$ (onto $\mc{N}_j$) respectively. It is easy to show that the restriction of $\psi_j$ to $\mc{M}_i$ is a normal endomorphism. We can write $\mc{M}_i$ in terms of the composite maps $\pi_i\psi_j$ in two ways:
\begin{align}
\mc{M}_i &= \pi_i(\mc{L}) = \pi_i\big(\psi_1(\mc{L}) \oplus \cdots \oplus \psi_s(\mc{L})\big) = \pi_i\psi_1(\mc{L}) + \cdots + \pi_i\psi_s(\mc{L}),\label{mdec1}\\
\mc{M}_i &= \pi_i(\mc{M}_i) = \pi_i\big(\psi_1(\mc{M}_i) \oplus \cdots \oplus \psi_s(\mc{M}_i)\big) = \pi_i\psi_1(\mc{M}_i) + \cdots + \pi_i\psi_s(\mc{M}_i).\label{mdec2}
\end{align}
(Here the vector space sum $+$ is not necessarily a direct sum.) We use $\pi_i\psi_j$ to denote both this map on $\mc{L}$ and its restriction to $\mc{M}_i$, as the meaning will be clear from the context; each composite projection $\pi_i\psi_j$ is a normal endomorphism of both $\mc{L}$ and $\mc{M}_i$.

\begin{lemma}\label{L2}
For each $i \in \{1, \ldots, r\}$ there exists $j \in \{1, \ldots, s\}$ such that:
\[
\pi_i(\mc{N}_j) = \mc{M}_i, \quad\ \mc{N}_j=\psi_j(\mc{M}_i), \quad\ \mc{N}_j \cong \mc{M}_i, \quad\ {\mc{N}_j}\!' = \mc{M}_i\!', \quad\ \mc{N}_j \s \mc{M}_i \oplus \mc{Z}(\,\overline{\mc{M}}_i).
\]
If $\mc{M}_i$ is non-abelian, $j$ is uniquely determined by $i$.
\end{lemma}

\noindent{\bf Proof}\, 
For each $n \in \N, \, (\pi_i\psi_j)^n = \pi_i(\psi_j\pi_i)^{n-1}\psi_j$ and $(\psi_j\pi_i)^n = \psi_j(\pi_i\psi_j)^{n-1}\pi_i$. Therefore, by Lemma \ref{L1}, the normal endomorphisms $\pi_i\psi_j : \mc{M}_i \to \mc{M}_i$ and $\psi_j\pi_i : \mc{N}_j \to \mc{N}_j$ are either both nilpotent or both bijective.  

If $\mc{M}_i$ is abelian then $\text{dim}(\mc{M}_i) = 1$.  From (\ref{mdec2}), there exists $j$ such that $\pi_i\psi_j(\mc{M}_i) \ne \{0\}$. For this $j$, therefore, $\pi_i\psi_j(\mc{M}_i) = \mc{M}_i$, and so $\pi_i\psi_j$ is bijective on $\mc{M}_i$. 

Now assume that $\mc{M}_i$ is non-abelian. Suppose that $\pi_i\psi_j$ is nilpotent on $\mc{M}_i$. Then $\iota - \pi_i\psi_j$ is bijective on $\mc{M}_i$, and therefore
\[
\mc{M}_i = (\iota - \pi_i\psi_j)(\mc{M}_i) = \pi_i(\iota - \psi_j)(\mc{M}_i)\s \pi_i(\,\overline{\mc{N}}_{\!j}).
\]
Thus
\[
\mc{M}_i\!' \s [\pi_i(\,\overline{\mc{N}}_{\!j}), \mc{M}_i] = [\overline{\mc{N}}_{\!j}, \mc{M}_i] \s \overline{\mc{N}}_{\!j}
\]
and hence $\psi_j(\mc{M}_i\!') = \{0\}$.  However,
\[
\psi_1(\mc{M}_i\!') + \cdots + \psi_S(\mc{M}_i\!') = \mc{M}_i\!' \ne \{0\},
\]
so there exists $j$ such that $\pi_i\psi_j$ is bijective, rather than nilpotent, on $\mc{M}_i$.

We have established that, whether or not $\mc{M}_i$ is abelian, there exists $j$ such that $\pi_i\psi_j$ and $\psi_j\pi_i$ are bijective;  the restricted projections for this $j$, namely
\[
\pi_i:\mc{N}_j\rightarrow\mc{M}_i,\qquad \psi_j:\mc{M}_i\rightarrow\mc{N}_j,
\]
are isomorphisms. Thus
\[
{\mc{N}_j}\!' = [\mc{N}_j,\psi_j(\mc{M}_i)] = \psi_j\big([\mc{N}_j, \mc{M}_i]\big) = [\mc{N}_j, \mc{M}_i] \s \mc{M}_i\!',
\]
and similarly $\mc{M}_i\!' \s {\mc{N}_j}\!'$, so $\mc{N}_j\!' = {\mc{M}_i}\!'$.  If $\mc{M}_i$ is non-abelian then ${\mc{M}_i}\!'\neq \{0\}$, in which case $j$ is uniquely determined by $i$.  Finally,
\beq\label{njz}
[(\iota - \pi_i)(\mc{N}_j), \overline{\mc{M}}_i] = [(\iota - \pi_i)(\mc{N}_j), \mc{L}] = (\iota - \pi_i)[\mc{N}_j, \mc{L}] = (\iota - \pi_i)({\mc{N}_j}\!') = (\iota - \pi_i)(\mc{M}_i\!') = \{0\},
\eeq
so 
$(\iota - \pi_i)(\mc{N}_j) \s \mc{Z}(\,\overline{\mc{M}}_i)$ and hence $\mc{N}_j \s \mc{M}_i \oplus \mc{Z}(\,\overline{\mc{M}}_i)$.  \hfill $\Box$ 

\begin{theorem}[Krull-Schmidt Theorem for Lie algebras]\label{T1}
~

\noindent With the above notation, $s = r$ and the summands can be numbered such that, for $i = 1, \ldots, r$, 
\beq\label{ks1}
\pi_i(\mc{N}_i) = \mc{M}_i, \quad\  \mc{N}_i=\psi_i(\mc{M}_i), \quad\ \mc{N}_i \cong \mc{M}_i, \quad\ {\mc{N}_i}' = \mc{M}_i\!',\quad\ \mc{N}_i \s \mc{M}_i \oplus \mc{Z}(\,\overline{\mc{M}}_i).
\eeq 
Furthermore, for each $k$ from $1$ to $r-1$,
\beq\label{ks2}
\mc{L} = \mc{M}_1 \oplus \cdots \oplus \mc{M}_k \oplus \mc{N}_{k+1} \oplus \cdots \oplus \mc{N}_r.
\eeq
If $\mc{Z}(\mc{L}) = \{0\}$ or $\mc{L}' = \mc{L}$, the decomposition (\ref{dsd}) is unique up to the order of the summands.
\end{theorem}  

\noindent{\bf Proof}\, By Lemma \ref{L2}, there exists $j$ such that $\pi_1\psi_j$ is bijective on $\mc{M}_1$ (and $\psi_j\pi_1$ is bijective on $\mc{N}_j$); if $\mc{M}_1$ is abelian, pick any such $j$. Choose the labels so that $\mc{N}_j$ becomes $\mc{N}_1$. 
Then
\[
(\pi_1\psi_1)^2(\mc{L}) = \pi_1\psi_1(\mc{M}_1) = \mc{M}_1,\quad\ \text{Ker}\big((\pi_1\psi_1)^2\big) = \text{Ker}(\pi_1\psi_1) = \text{Ker}(\psi_1\pi_1\psi_1) = \text{Ker}(\psi_1) = \overline{\mc{N}}_{\!1},
\]
 so, by Fitting's Lemma, $\mc{L} = \mc{M}_1 \oplus \overline{\mc{N}}_{\!1}$.  Hence (\ref{ks2}) holds when $k = 1$.
Now suppose that
\[
\mc{L} = \mc{M}_1 \oplus \cdots \oplus \mc{M}_k \oplus \mc{N}_{k+1} \oplus \cdots \oplus \mc{N}_r
\]
for a particular $k<r$. Comparing this decomposition with 
\[
\mc{L} = \mc{M}_1 \oplus \cdots \oplus \mc{M}_k \oplus \mc{M}_{k+1} \oplus \cdots \oplus \mc{M}_r,
\]
the above reasoning enables us to identify a summand in the former decomposition, having all of the stated properties, which can be replaced by $\mc{M}_{k+1}$.  This is not any of $\mc{M}_1, \ldots, \mc{M}_k$ because $\pi_{k+1}(\mc{M}_i) = 0$ for $i \ne k+1$; so, renumbering if necessary, we can take it to be $\mc{N}_{k+1}$.  Then
\[
\mc{L} = \mc{M}_1 \oplus \cdots \oplus \mc{M}_k \oplus \mc{M}_{k+1} \oplus \mc{N}_{k+2} \oplus \cdots \oplus \mc{N}_r.
\]
By induction, (\ref{ks2}) holds for $k = 2, \ldots, r$ and both decompositions contain the same number of summands.
At each stage, we have arranged the remaining labels so that $j=i$ and hence (\ref{ks1}) holds. If $\mc{Z}(\mc{L}) = \{0\}$ then $\mc{M}_i = \mc{N}_i$ for each $i$.  If $\mc{L}' = \mc{L}$ then $\mc{M}_i = \mc{M}_i\!' = \mc{N}_i\!\,' = \mc{N}_i$.  In each of these cases, $\mc{L}$ has a unique decomposition (up to the ordering of the summands).  \hfill $\Box$

\begin{theorem}\label{T2}
A bijective linear map $\phi : \mc{L} \to \mc{L}$ is an automorphism of $\mc{L}$ if and only if it has the form
$\theta + \zeta$, where $\theta$ is an automorphism of $\mc{L}$ that maps each $\mc{M}_i$ to itself or to an isomorphic summand, $\mc{M}_j$, and where $\zeta$ is a linear map of $\mc{L}$ such that $\zeta(\mc{L}) \s \mc{Z}(\mc{L})$ and $\zeta(\mc{L}') = \{0\}$. 
\end{theorem}

\noindent{\bf Proof}\, Any automorphism $\phi$ of (\ref{dsd}) gives a decomposition $\mc{L} = \phi(\mc{M}_1) \oplus \cdots \oplus \phi(\mc{M}_r)$. Label $\phi(\mc{M}_1), \ldots, \phi(\mc{M}_r)$ as $\mc{N}_1, \ldots, \mc{N}_r$ in such an order that $\mc{N}_i$ corresponds to $\mc{M}_i$ as in Theorem \ref{T1}; then there is a permutation $p(i)$ of $(1,\dots,r)$ such that $\mc{N}_i=\phi(\mc{M}_{p(i)})$ and $\mc{M}_{p(i)}\cong \mc{M}_i$.
Now define the linear maps $\zeta$ and $\theta$ on $\mc{L}$ by
\[
\theta(X) = \sum\limits_{i=1}^r \pi_i\phi\:\!\pi_{p(i)}(X),\qquad \zeta(X) = \sum\limits_{i=1}^r (\iota - \pi_i)\:\!\phi\:\!\pi_{p(i)}(X),\qquad X \in \mc{L}.
\]
Clearly, $\phi = \theta + \zeta$; morever, $\zeta(\mc{L}) \s \mc{Z}(\mc{L})$ because, from (\ref{njz}), 
\[
(\iota - \pi_i)\:\!\phi\:\!\pi_{p(i)}(X) \in (\iota - \pi_i)\:\!\mc{N}_i \s \mc{Z}(\mc{L}),\qquad X \in \mc{L}.
\]
Note that $\phi(\mc{M}_{p(i)}\!\!\!\!\!\!\!'\ \,\;)=\mc{N}_i\,\!'=\mc{M}_i\!'$, and so $\zeta(\mc{L}') = \{0\}$. By definition,
\[
\theta\big(\mc{M}_{p(i)}\big)=\pi_i\phi\big(\mc{M}_{p(i)}\big)=\pi_i\,\mc{N}_i=\mc{M}_i,
\]
so $\theta(\mc{L}) = \mc{L}$; as $\theta$ is a composite of Lie algebra homomorphisms, it is an automorphism of $\mc{L}$. 

Conversely, suppose that $\theta$ and $\zeta$ have the stated properties. Then $\theta$ permutes the components $\mc{M}_i$; it is of the form $\theta_1 + \cdots + \theta_r$, where each $\theta_i$ is an isomorphism from $\mc{M}_i$ to some $\mc{M}_j$.  
As $\phi = \theta + \zeta$ is the sum of two linear maps of $\mc{L}$, it is a linear map of $\mc{L}$.   For any $X, Y \in \mc{L}$,
\[
\phi([X,Y])\! =\! \theta([X,Y]) + \zeta([X,Y])\! =\! [\theta(X),\theta(Y)]\! =\! [\theta(X) + \zeta(X), \theta(Y) + \zeta(Y)]\! =\! [\phi(X),\phi(Y)].
\]
Thus $\phi$ is an endomorphism of $\mc{L}$; provided that $\phi$ is bijective, it is an automorphism of $\mc{L}$. \hfill $\Box$ \\

When $\mc{L}$ is finite-dimensional, Theorem \ref{T2} gives necessary and sufficient conditions on the matrix representing an automorphism of $\mc{L}$.
Let $B = (b_{p}^q)$ be the matrix of $\phi$ relative to a basis for $\mc{L}$ ordered in such a way that a basis for $\mc{M}_1$ is followed by a basis for $\mc{M}_2$, etc.  
Let $\{X_{i_1}, \ldots, X_{i_{k_i}}\}$ be a basis for $\mc{M}_i$.
If it is not clear how $\mc{L}$ can be expressed as a direct sum, the algorithm in \cite{deGraaf} can be used. Note that this may reveal the presence of a central component, in which case one or more of the direct summands will be one-dimensional.  Their presence does not affect our argument.

If $\mc{M}_i$ is not isomorphic to any other direct summand then the sub-matrix of $B$ consisting of the $b_{p}^q$ for
$i_1 \le p \le  i_{k_i}, i_1 \le q \le  i_{k_i}$ represents an automorphism of $\mc{M}_i$.
All entries are zero outside the blocks corresponding to automorphisms of the $\mc{M}_i$, except when 
$X_p \notin \mc{L}', X_q \in \mc{Z}(\mc{L})$, in which case $b_{p}^q$ is arbitrary subject to $B$ being non-singular.

If $\mc{M}_i \cong \mc{M}_{i'}$ for some $i' \ne i$, the above automorphisms may be composed with an automorphism that swaps $\mc{M}_i$ and $\mc{M}_{i'}$. So whenever an irreducible decomposition of $\mc{L}$ has two or more isomorphic summands, any permutation of these is allowable.

{\bf Example}\, Let $\mc{L}$ have a basis $\{X_1, \cdots, X_8\}$ with non-zero commutator relations
\[[X_1,X_2]=X_1, \quad [X_4,X_5]=X_3, \quad [X_7,X_8]=X_6.
\]
Then $\mc{L}$ is the direct sum $\mc{M}_1 \oplus \mc{M}_2 \oplus \mc{M}_3$ where
$\mc{M}_1 \cong A_{2,1}$, $\mc{M}_2 \cong \mc{M}_3 \cong A_{3,1}$ (see Table 1). 
The automorphisms of $A_{2,1}$, with $[X_1,X_2]=X_1$, have matrices of the form
\[
\left(\begin{array}{cc} a & 0 \\ b & 1 \end{array}\right); \quad a\ne 0, \quad a,b\in \R.
\]
The automorphisms of the nilpotent Lie algebra $A_{3,1}$, with $[X_4,X_5]=X_3$, have matrices of the form 
\[
\left(\begin{array}{ccc}eh - fg & 0 & 0 \\ c & e & g \\ d & f & h \end{array}\right); \quad eh-fg \ne 0, \quad c,d,e,f,g,h\in\R.
\]
Moreover, $\mc{Z}(\mc{L})$ = span$\{X_3, X_6\}$ and $[\mc{L},\mc{L}]$ = span$\{X_1,X_3,X_6\}$.
Hence any matrix of the following form represents
an automorphism of $\mc{L}$:
\[
\left(\begin{array}{cccccccc} a & 0 & 0 & 0 & 0 & 0 & 0 & 0 \\ b & 1 & \alpha & 0 & 0 & \beta & 0 & 0 \\
0 & 0 & eh - fg & 0 & 0 & 0 & 0 & 0 \\ 0 & 0 & c & e & g & \gamma & 0 & 0 \\ 0 & 0 & d & f & h & \delta & 0 & 0 \\
0 & 0 & 0 & 0 & 0 & kn - lm & 0 & 0 \\ 0 & 0 & \varepsilon & 0 & 0 & i & k & m \\
0 & 0 & \zeta & 0 & 0 & j & l & n \end{array}\right),\qquad a(eh-fg)(kn-lm)\ne 0. 
\]
Here Greek letters represent the additional arbitrary $\mc{Z}(\mc{L})$ components.
As $\mc{M}_2 \cong \mc{M}_3$, there are other automorphisms formed by swapping these two direct summands; their matrices are
\[
\left(\begin{array}{cccccccc} a & 0 & 0 & 0 & 0 & 0 & 0 & 0 \\ b & 1 & \alpha & 0 & 0 & \beta & 0 & 0 \\
0 & 0 & 0 & 0 & 0 & eh - fg & 0 & 0 \\ 0 & 0 & \gamma & 0 & 0 & c & e & g  \\ 0 & 0 & \delta & 0 & 0& d & f & h \\
0 & 0 & kn - lm & 0 & 0 & 0 & 0 & 0 \\ 0 & 0 & i & k & m & \varepsilon & 0 & 0 \\
0 & 0 & j & l & n & \zeta & 0 & 0 \end{array}\right),\qquad a(eh-fg)(kn-lm)\ne 0. 
\]
By Theorem \ref{T2}, there are no other automorphisms of $\mc{L}$.

\section{Automorphisms of real Lie algebras of dimension $\leq 5$}

We now present the classification of automorphism groups for indecomposable real Lie algebras of dimension five or less. The Lie algebras are labelled according to the tables in the classification of invariants by Patera \textit{et al.} \cite{PSWH}. In a few instances, it is helpful to use a basis that is slightly different to \cite{PSWH}; where this occurs, we append an asterisk, $(*)$, to the label.

\subsection{How to use the tables}

Tables 1, 2 and 3 are designed to present the results concisely, avoiding the need to write out matrices explicitly (except in one particularly complicated case which is treated separately). Here is a typical entry:

{\small
\[
\begin{array}{lllcccc}
\mathrm{Name}\quad&\mathrm{Notes}\quad&\mathrm{Nonzero}\ c_{ij}^k \ (i<j)\quad&&\mathrm{Discrete\ Generators}\quad&\mathrm{Extra\ Outer\ Der.}&\mathrm{Block\ Diag.}\\[2.5ex]
A_{4,8}&&c_{23}^1\!=\!c_{24}^2\!=\!1,\ c_{34}^3\!=\!-1&&p_{12},\ (-X_1,X_3,X_2,-X_4)&E_1^1\!+\!E_3^3,\ E_4^1&(1,1,1,1)\\[0.9ex]
\end{array}
\]}

The first three columns give the name of the Lie algebra (in this case, $A_{4,8}$), any explanatory notes (such as restrictions on parameter values), and all nonzero structure constants $c_{ij}^k$ for which $i<j$. To understand the remaining columns, it is helpful to look in detail at the automorphisms of $A_{4,8}$. For this Lie algebra, the determining equations \eqref{AC} are easily solved without the aid of inner automorphisms. There are two families of solutions:
\beq\label{autex}
B_1 = \left[
\begin{array}{cccc}
b^2_2b_3^3 & 0 & 0 & 0\\
b^2_2b^3_4 & b_2^2 & 0 & 0\\
b_3^3b_4^2 & 0 & b_3^3 &0\\
b_4^1 & b_4^2 & b_4^3 & 1
\end{array}
\right], \ b_2^2b_3^3\neq 0; \qquad B_2 = \left[
\begin{array}{cccc}
-b_2^3b_3^2 & 0 & 0 & 0\\
-b_2^3b^2_4 & 0 & b_2^3 & 0\\
-b_3^2b_4^3 & b_3^2 & 0 &0\\
b_4^1 & b_4^2 & b_4^3 & -1
\end{array}
\right],\ b_2^3b_3^2\neq 0;
\eeq
all parameters are arbitrary real numbers, subject only to the determinant of each matrix being nonzero. Although we have not needed to use the inner automorphisms so far, it is useful to do so now. As $X_1$ belongs to the centre of 
$A_{4,8}$, the matrix $A_1(\varepsilon_1)$ is the identity matrix and cannot change $B$. The remaining matrices that generate the inner automorphisms are
\[
A_2(\varepsilon_2) = \!\left[
\begin{array}{cccc}
1 & 0 & 0 & 0\\
0 & 1 & 0 & 0\\
-\varepsilon_2 & 0 & 1 &0\\
0 & -\varepsilon_2 & 0 & 1
\end{array}
\right]\!\!, \ \,
A_3(\varepsilon_3) = \!\left[
\begin{array}{cccc}
1 & 0 & 0 & 0\\
\varepsilon_3 & 1 & 0 & 0\\
0 & 0 & 1 &0\\
0 & 0 & \varepsilon_3 & 1
\end{array}
\right]\!\!, \ \,
A_4(\varepsilon_4) = \!\left[
\begin{array}{cccc}
1 & 0 & 0 & 0\\
0 & e^{\varepsilon_4} & 0 & 0\\
0 & 0 & e^{-\varepsilon_4} &0\\
0 & 0 & 0 & 1
\end{array}
\right]\!\!,
\]
where each $\varepsilon_j$ is a real-valued parameter. To simplify the automorphism matrices \eqref{autex}, let
\[
\tilde{B}_1=A_4(-\ln|b_2^2|)A_2(b_3^1/b_1^1)A_3(-b_2^1/b_1^1)B_1,\qquad
\tilde{B}_2=A_4(\ln|b_3^2|)A_2(b_3^1/b_1^1)A_3(-b_2^1/b_1^1)B_2.
\]
After a change of parameters to remove unnecessary indices, we obtain
\beq\label{autexout}
\tilde{B}_1 = \left[
\begin{array}{cccc}
 \epsilon a & 0 & 0 & 0\\
0 & \epsilon & 0 & 0\\
0 & 0 & a &0\\
b & 0 & 0 & 1
\end{array}
\right],\qquad \tilde{B}_2 = \left[
\begin{array}{cccc}
 -\epsilon a & 0 & 0 & 0\\
0 & 0 & a & 0\\
0 & \epsilon & 0 & 0\\
-b & 0 & 0 & -1
\end{array}
\right],\qquad \epsilon=\pm 1,\ a\neq 0.
\eeq
All that remains is to decompose these automorphisms, as follows:
\beq\label{tbdec}
\tilde{B}=\Delta\,\exp\left[
\begin{array}{cccc}
0 & 0 & 0 & 0\\
0 & 0 & 0 & 0\\
0 & 0 & 0 &0\\
\beta & 0 & 0 & 0
\end{array}
\right]\exp
\left[
\begin{array}{cccc}
\alpha & 0 & 0 & 0\\
0 & 0 & 0 & 0\\
0 & 0 & \alpha & 0\\
0 & 0 & 0 & 0
\end{array}
\right],
\eeq
where the matrix $\Delta$ is an element of the discrete subgroup generated by
\beq\label{discex}
\Delta_1 = \left[
\begin{array}{cccc}
-1 & 0 & 0 & 0\\
0 & -1 & 0 & 0\\
0 & 0 & 1 & 0\\
0 & 0 & 0 & 1
\end{array}
\right],\qquad
\Delta_2 = \left[
\begin{array}{cccc}
-1 & 0 & 0 & 0\\
0 & 0 & 1 & 0\\
0 & 1 & 0 &0\\
0 & 0 & 0 & -1
\end{array}
\right].
\eeq
This eight-element subgroup is isomorphic to the dihedral group $D_4$.  At this stage, we have decomposed the automorphism group into elements that can be described with fairly simple notation, as follows.
\begin{enumerate}
\item Let $p_\mathbf{m}$ (where $\mathbf{m}$ is a multi-index) denote the matrix $\mathrm{diag}\{s_1,\dots,s_R\}$, where $s_i=-1$ if $i$ is in $\mathbf{m}$ and $s_i=1$ otherwise.  For $R<10$, the components of the multi-index are single digits and, to save space, we will not separate them by commas. So in our example, $\Delta_1=p_{12}$.
\item Denote a discrete symmetry that is not of the above form,
\[
\Phi:(X_1,\dots,X_R)\longmapsto(\tilde{X}_1,\dots,\tilde{X}_R)=(b_1^lX_l,\dots,b_R^lX_l),
\]
by $(\tilde{X}_1,\dots,\tilde{X}_R)$ expressed in terms of the generators $X_i$. In our example,
\[
\Delta_2=(-X_1,X_3,X_2,-X_4).
\]
\item Where there is a one-parameter (local) subgroup of exponentiated outer derivations, write the generator in terms of the Weyl basis; here $E_i^j$ is the matrix whose only nonzero entry is a $1$ in row $i$, column $j$. In our example, the subgroup that is parametrized by $\alpha$ in \eqref{tbdec} is generated by $E_1^1+ E_3^3$, while the subgroup that is parametrized by $\beta$ is generated by $E_4^1$.
\end{enumerate}

In the table, generators of a discrete group of automorphisms are listed in the column headed `Discrete Generators', and generators of one-parameter local Lie subgroups of exponentiated outer derivations are placed under the heading `Extra Outer Derivations'. In principle, one can list all components of $\tilde{B}$ under one of these two headings (omitting the word `Extra'). However, this turns out to be unwieldy for Lie algebras with large automorphism groups. As a compromise, we use a sixth column, as follows.

\begin{enumerate}\setcounter{enumi}{3}
\item The `Block Diagonal' entry lists the block diagonal structure, ordered by row/column number; for Lie algebras other than $A_{5,4}$ (see below), all entries are nonzero scalars and/or $\mathrm{SL}(n,\mathbb{R})$ matrix blocks. Each scalar is either $1$ or an arbitrary nonzero real number (denoted $a,b,c,\dots$); the identity matrix corresponds to a string of $1$s. Where an $\mathrm{SL}(n,\mathbb{R})$ block is entirely arbitrary (apart from the condition that its determinant is $1$), it is denoted by $S_\mathbf{m}$, where the multi-index $\mathbf{m}$ lists the $n$ contiguous rows/columns where the block is located. If the same block occurs elsewhere in the block diagonal structure, the multi-index will be repeated.
For instance, the block diagonal matrices
\[
\left[
\begin{array}{cccc}
a^2 & 0 & 0 & 0\\
0 & b & 0 & 0\\
0 & 0 & ab & 0\\
0 & 0 & 0 & c
\end{array}
\right],\quad\ 
\left[
\begin{array}{cccc}
1 & 0 & 0 & 0\\
0 & p & q & 0\\
0 & r & s & 0\\
0 & 0 & 0 & 1
\end{array}
\right],\quad\ 
\left[
\begin{array}{cccc}
p & q & 0 & 0\\
r & s & 0 & 0\\
0 & 0 & ap & aq\\
0 & 0 & ar & as
\end{array}
\right],\quad\ \text{where}\ ps-qr=1,
\]
are denoted respectively by
\[
(a^2,b,ab,c),\qquad (1,S_{23},1),\qquad (S_{12},aS_{12}).
\]
For $A_{5,4}$, the $4\times 4$ block labelled $\mathrm{Sp}(4,\mathbb{R})^T$ denotes the \textit{transpose} of an arbitrary matrix $M$ in the symplectic group $\mathrm{Sp}(4,\mathbb{R})$. So $M$ is subject only to the constraint
\[
M^T\Omega M=\Omega,\quad\text{where}\quad\Omega =\left[
\begin{array}{cccc}
0 & 0 & 1 & 0\\
0 & 0 & 0 & 1\\
-1 & 0 & 0 & 0\\
0 & -1 & 0 & 0
\end{array}
\right].
\]

\end{enumerate} 
The advantage of the above compromise is that the notation remains concise. One can decompose the block diagonal structure, but this is not necessary for the reconstruction of the full automorphism group, which proceeds as follows.

Premultiply the block diagonal matrix by each exponentiated extra outer derivation. Premultiply the result by an element $\Delta$ of the discrete group whose generators are in the fourth column; one needs to keep track of the separate outcomes of using each possible $\Delta$. Finally, premultiply every outcome by the inner automorphisms $A_1(\varepsilon_1)\cdots A_R(\varepsilon_R)$, where each $\varepsilon_j$ takes all possible values. Schematically, each $\Delta$ produces a matrix of automorphisms,
\[
B = \big( \prod A_j(\varepsilon_j)  \big)\ \Delta\ \big(\prod\text{exp}(\alpha_i\,\text{Extra Outer Derivation}_i)  \big) \big(\text{Block Diagonal} \big);
\]
the set of all such matrices comprises the full automorphism group.

The above information is sufficient to construct the automorphism group of each indecomposable Lie algebra of dimension $\leq 5$ (bar one) from the tables that follow (the sole exception is treated separately in \S4.2). However, we have added two further pieces of notation, for clarity.
 
\begin{enumerate}\setcounter{enumi}{4}
\item The `Discrete Generators' lists for $A_{3,8}$ and $A_{5,40}$ each include a Weyl reflection, which is denoted by double parentheses. This is an inner automorphism that cannot be expressed as a product in which each $A_j$ is used once only. Weyl reflections are the only inner automorphisms that appear in the last three columns.

\item Where an extra outer derivation is enclosed in square brackets with a subscript $u$ (for example, $[E_1^1+E_2^2]_u$), one can exponentiate it with a parameter that is restricted to be in the range $[0,2\pi |u|)$. The remaining values are generated by inner automorphisms. However, if one is only interested in reconstructing the full automorphism group, it is simpler to allow the parameter to take arbitrary values, while recognising that the inner automorphisms will then produce some duplications.
\end{enumerate}

\subsection{The automorphism group of $A_{5,17}^{u,v,w}$}
The nonzero structure constants $c_{ij}^k$ $(i<j)$ are
\[
c_{15}^1=c_{25}^2=u,\;c_{35}^3=c_{45}^4=v,\;c_{15}^2=-1,\;c_{25}^1=1,\;c_{35}^4=-w,\;c_{45}^3=w,\qquad w\neq 0.
\]
For all parameter values, every real nonsingular matrix of the form
\[B_1=\left[
\begin{array}{ccccc}
a&b&0&0&0\\
-b&a&0&0&0\\
0&0&g&h&0\\
0&0&-h&g&0\\
k_1&k_2&k_3&k_4&1
\end{array}
\right]
\]
is an automorphism; indeed,  every automorphism matrix is of this form, except in the following special cases.

If $u=-v\neq 0$ and $\left|w\right|=1 $ then the group of automorphisms is generated by $B_1$ and
\[B_2=\left[
\begin{array}{ccccc}
0&0&0&w&0\\
0&0&1&0&0\\
0&w&0&0&0\\
1&0&0&0&0\\
0&0&0&0&-1
\end{array}
\right].
\]

If $u=v\neq 0$ and $\left|w\right|=1$ then the group of automorphisms consists of all nonsingular matrices of the following form:
\[B_3=\left[
\begin{array}{ccccc}
a&b&c&d&0\\
-b&a&-wd&wc&0\\
e&f&g&h&0\\
-wf&we&-h&g&0\\
k_1&k_2&k_3&k_4&1
\end{array}
\right],
\]

If $u=v=0$ there is an extra subgroup of automorphisms that is generated by $p_{245}$, so the automorphism group is generated by: (a) $B_1$ and $p_{245}$, if $|w|\neq 1$; (b) $B_3$ and $p_{245}$, if $|w|=1$.

\renewcommand\arraystretch{1.4} 

\newgeometry{bottom=3cm}

\setcounter{LTchunksize}{100}
\begin{landscape}
\begin{center}

\scriptsize
\rowcolors{1}{light-gray}{white}
\begin{longtable}{lllccc}
\caption{\bfseries Automorphisms of $2$- and $3$-dimensional indecomposable real Lie algebras}\\
Name\quad&Notes\qquad&Structure Constants $c_{ij}^k \
(i<j)$\quad&Discrete\ Generators\quad&Extra\ Outer\ Derivations\quad&Block\ Diagonal\\[0.1cm] \hline\vspace*{-0.38cm}
\endfirsthead
Name\quad&Notes\qquad&Structure Constants $c_{ij}^k \
(i<j)$\quad&Discrete\ Generators\quad&Extra\ Outer\ Derivations\quad&Block\ Diagonal\\[0.1cm] \hline\vspace*{-3.4ex}
\endhead
\hline \multicolumn{6}{r}{\emph{Continued on next page}}
\endfoot
\hline
\endlastfoot
$A_{2,1}$   &  & $c_{12}^1=1$ & $p_1$&   & $(1,1)$ \\[0.9ex]
$A_{3,1}$   & $\mathrm{nilpotent}$  & $c_{23}^1=1$  &$p_{12}$& $E_1^1\!+\!E_2^2$  & $(1,S_{23})$ \\[0.9ex]
$A_{3,2}$   & & $c_{13}^1=c_{23}^1=c_{23}^2=1$    && & $(a,a,1)$ \\[0.9ex]
$A_{3,3}$   &  & $c_{13}^1=c_{23}^2=1$  & $p_1$&   & $(S_{12},1)$ \\[0.9ex]
$A_{3,4}$   &  & $c_{13}^1=1,c_{23}^2=-1$   & $(-\!X_2,X_1,-\!X_3)$&& $(1,a,1)$ \\[0.9ex]
$A_{3,5}^u$ & $0\!<\!|u|\!<\!1$  & $c_{13}^1=1,c_{23}^2=u$      & $p_1$&  & $(1,a,1)$ \\[0.9ex]
$A_{3,6}$   &  & $c_{13}^2=-1,c_{23}^1=1$  & $p_{23}$& $E_1^1\!+\!E_2^2$ & $(1,1,1)$ \\[0.9ex]
$A_{3,7}^u$ & $u\!>\!0$  & $c_{13}^1=c_{23}^2=u,c_{13}^2=-1,c_{23}^1=1$ && $\big[ E_1^1\!+\!E_2^2\big]_u$ & $(1,1,1)$ \\[0.9ex]
$A_{3,8}$   & $\mathfrak{sl}(2,\mathbb{R})$ & $c_{12}^1=c_{23}^3=1,c_{13}^2=-2$ & $p_{13},\ ((-\!X_3,-\!X_2,-\!X_1))$& & $(1,1,1)$ \\[0.9ex]
$A_{3,9}$   & $\mathfrak{so}(3)$            & $c_{12}^3=c_{23}^1=1,c_{13}^2=-1$ &&                           & $(1,1,1)$ \\[0.9ex]
\end{longtable}
\end{center}
\normalsize

\setcounter{LTchunksize}{100}
\begin{center}
\scriptsize
\rowcolors{1}{light-gray}{white}
\begin{longtable}{llllccc}
\caption{\bfseries Automorphisms of $4$-dimensional indecomposable real Lie algebras}\\
Name\quad&Notes\qquad&Structure Constants $c_{ij}^k \
(i<j)$\quad&$\phantom{0}$&Discrete\ Generators\quad&Extra\ Outer\ Derivations\quad&Block\ Diagonal\\[0.1cm] \hline\vspace*{-0.38cm}
\endfirsthead
Name\quad&Notes\qquad&Structure Constants $c_{ij}^k \
(i<j)$\quad&$\phantom{0}$&Discrete\ Generators\quad&Extra\ Outer\ Derivations\quad&Block\ Diagonal\\[0.1cm] \hline\vspace*{-3.4ex}
\endhead
\hline \multicolumn{7}{r}{\emph{Continued on next page}}
\endfoot
\hline
\endlastfoot
$A_{4,1}$&$\text{nilpotent}$&$c_{24}^1\!=\!c_{34}^2\!=\!1$&&&$E_3^1,\ E_4^3$&$(ab^2,ab,a,b)$\\[0.9ex]
$A_{4,2}^u$&$u\!\not\in\! \{0,1\}$&$c_{14}^1\!=\!u,\ c_{24}^2\!=\!c_{34}^2\!=\!c_{34}^3\!=\!1$&&&&$(a,b,b,1)$\\[0.9ex]
$A_{4,2}^1$&&$c_{14}^1\!=\!c_{24}^2\!=\!c_{34}^2\!=\!c_{34}^3\!=\!1$&&&$E_1^2,\ E_3^1$&$(a,b,b,1)$\\[0.9ex]
$A_{4,3}$&&$c_{14}^1\!=\!c_{34}^2\!=\!1$&&&$E_4^3$&$(a,b,b,1)$\\[0.9ex]
$A_{4,4}$&&$c_{14}^1\!=\!c_{24}^1\!=\!c_{24}^2\!=\!c_{34}^2\!=\!c_{34}^3\!=\!1$&&&$E_3^1$&$(a,a,a,1)$\\[0.9ex]
$A_{4,5}^{u,v}$&$uv\!\neq\! 0,\ -1\!\leq\! u\!<\!v\!<\!1$&$c_{14}^1\!=\!1,\ c_{24}^2\!=\!u,\ c_{34}^3\!=\!v$&&$p_1$&&$(1,a,b,1)$\\[0.9ex]
$A_{4,5}^{u,u}$&$u\!\neq\! 0,\ -1\!\leq\! u\!<\!1$&$c_{14}^1\!=\!1,\ c_{24}^2\!=\!c_{34}^3\!=\!u$&&$p_1,p_2$&$E_2^2$&$(1,S_{23},1)$\\[0.9ex]
$A_{4,5}^{u,1}(*)$&$u\!\neq\! 0,\ -1\!\leq\! u\!<\!1$&$c_{14}^1\!=\!u,\ c_{24}^2\!=\!c_{34}^3\!=\!1$&&$p_1,p_2$&$E_2^2$&$(1,S_{23},1)$\\[0.9ex]
$A_{4,5}^{1,1}$&&$c_{14}^1\!=\!c_{24}^2\!=\!c_{34}^3\!=\!1$&&$p_1$&&$(S_{123},1)$\\[0.9ex]
$A_{4,6}^{u,v}$&$u\!\neq\! 0,\ v\!\geq\! 0$&$c_{14}^1\!=\!u,\ c_{24}^2\!=\!c_{34}^3\!=\!v,\ c_{24}^3\!=\!-1,\ c_{34}^2\!=\!1$&&&$[E_2^2\!+\!E_3^3]_v$&$(a,1,1,1)$\\[0.9ex]
$A_{4,7}$&&$c_{14}^1\!=\!2,\ c_{23}^1\!=\!c_{24}^2\!=\!c_{34}^2\!=\!c_{34}^3\!=\!1$&&&&$(a^2,a,a,1)$\\[0.9ex]
$A_{4,8}$&&$c_{23}^1\!=\!c_{24}^2\!=\!1,\ c_{34}^3\!=\!-1$&&$p_{12},\ (-X_1,X_3,X_2,-X_4)$&$E_1^1\!+\!E_3^3,\ E_4^1$&$(1,1,1,1)$\\[0.9ex]
$A_{4,9}^u$&$-1\!<\!u\!<\!1$&$c_{14}^1\!=\!u\!+\!1,\ c_{23}^1\!=\!c_{24}^2\!=\!1,\ c_{34}^3\!=\!u$&&$p_{12}$&&$(a,1,a,1)$\\[0.9ex]
$A_{4,9}^1$&&$c_{14}^1\!=\!2,\ c_{23}^1\!=\!c_{24}^2\!=\!c_{34}^3\!=\!1$&&$p_{12}$&&$(1,S_{23},1)$\\[0.9ex]
$A_{4,10}$&&$c_{23}^1\!=\!c_{34}^2\!=1,\ \!c_{24}^3\!=\!-1$&&$p_{124}$&$2E_1^1\!+\!E_2^2\!+\!E_3^3,\ E_4^1$&$(1,1,1,1)$\\[0.9ex]
$A_{4,11}^u$&$u\!>\!0$&$c_{14}^1\!=\!2u,\ c_{23}^1\!=\!c_{34}^2\!=\!1,\ c_{24}^2\!=\!c_{34}^3\!=\!u,\ c_{24}^3\!=\!-1$&&&$\big[ 2E_1^1\!+\!E_2^2\!+\!E_3^3 \big]_u$&$(1,1,1,1)$\\[0.9ex]
$A_{4,12}$&&$c_{13}^1\!=\!c_{23}^2\!=\!c_{24}^1\!=\!1,\ c_{14}^2\!=\!-1$&&$p_{24}$&&$(1,1,1,1)$
\end{longtable}
\end{center}
\normalsize

\setcounter{LTchunksize}{100}
\begin{center}
\scriptsize
\rowcolors{1}{light-gray}{white}
\begin{longtable}{llp{5.2cm}p{0cm}ccc}
\caption{\bfseries Automorphisms of $5$-dimensional indecomposable real Lie algebras}\\
Name\quad&Notes\qquad&Structure Constants $c_{ij}^k \
(i<j)$\quad&$\phantom{0}$&Discrete\ Generators\quad&Extra\ Outer\ Derivations\quad&Block\ Diagonal\\[0.1cm] \hline\vspace*{-0.38cm}
\endfirsthead
Name\quad&Notes\qquad&Structure Constants $c_{ij}^k \
(i<j)$\quad&$\phantom{0}$&Discrete\ Generators\quad&Extra\ Outer\ Derivations\quad&Block\ Diagonal\\[0.1cm] \hline\vspace*{-3.4ex}
\endhead
\hline \multicolumn{7}{r}{\emph{Continued on next page}}
\endfoot
\hline
\endlastfoot
$A_{5,1}$&\text{nilpotent}&$c_{35}^1\!=\!c_{45}^2\!=1$&&$p_{13}$&\begin{minipage}[c]{3.8cm}\begin{center}\vspace{1ex}$E_1^1\!+\!E_3^3,\ E_3^1,\ E_3^2$,\\[2pt] $E_4^1,\ E_4^2,\ E_5^3,\ E_5^4$\vspace{1ex}\end{center}\end{minipage}&$(aS_{34},S_{34},a)$\\[0.9ex]
$A_{5,2}$&\text{nilpotent}&$c_{25}^1\!=\!c_{35}^2\!=\!c_{45}^3\!=1$&&&$E_3^1\!+\!E_4^2,\ E_4^1,\ E_5^4$&$(ab^3,ab^2,ab,a,b)$\\[0.9ex]
$A_{5,3}(*)$&\text{nilpotent}&$c_{34}^1\!=\!c_{35}^2\!=\!c_{45}^3\!=1$&&$p_{234}$&\begin{minipage}[c]{3.8cm}\begin{center}\vspace{1ex}$2E_1^1\!+\!E_2^2\!+\!E_3^3\!+\!E_4^4$,\\[2pt] $E_4^1,\ E_4^2,\ E_5^1,\ E_5^2$\vspace{1ex}\end{center}\end{minipage}&$(S_{12},1,S_{12})$\\[0.9ex]
$A_{5,4}$&\text{nilpotent}&$c_{24}^1\!=\!c_{35}^1\!=1$&&$p_{123}$&$E_1^1\!+\!E_2^2\!+\!E_3^3$&$(1,\mathrm{Sp}(4,\mathbb{R})^T)$\\[0.9ex]
$A_{5,5}$&\text{nilpotent}&$c_{25}^1\!=\!c_{34}^1\!=\!c_{35}^2\!=1$&&&$E_3^2\!-\!E_5^4,\ E_3^4,\ E_4^2\!+\!E_5^3,\ E_5^2$&$(ab^2,ab,a,b^2,b)$\\[0.9ex]
$A_{5,6}$&\text{nilpotent}&$c_{25}^1\!=\!c_{34}^1\!=\!c_{35}^2\!=\!c_{45}^3\!=1$&&&$E_3^2\!+\!E_4^3\!-\!E_5^4,\ E_4^2\!+\!E_5^3,\ E_5^2$&$(a^5,a^4,a^3,a^2,a)$\\[0.9ex]
$A_{5,7}^{u,v,w}$ & $uvw\!\neq\!0,\ -\!1\!<\!w\!<\!v\!<\!u\!<\!1$ & $c_{15}^1\!=\!1,\ c_{25}^2\!=\!u,\ c_{35}^3\!=\!v,\ c_{45}^4\!=\!w$&&$p_1$&&$(1,a,b,c,1)$\\[0.9ex]
$A_{5,7}^{u,u,w}$&$uw\!\neq\!0,\ -\!1\!\leq\!w\!<\!u\!<\!1$&$c_{15}^1\!=\!1,\ c_{25}^2\!=\!c_{35}^3\!=\!u,\ c_{45}^4\!=\!w$&&$p_1,p_2$&$E_2^2$&$(1,S_{23},a,1)$\\[0.9ex]
$A_{5,7}^{u,v,v}$&$uv\!\neq\!0,\ -\!1\!\leq\!v\!<\!u\!<\!1$&$c_{15}^1\!=\!1,\ c_{25}^2\!=\!u,\ c_{35}^3\!=\!c_{45}^4\!=\!v$&&$p_1,p_3$&$E_3^3$&$(1,a,S_{34},1)$\\[0.9ex]
$A_{5,7}^{1,v,w}$&$vw\!\neq\!0,\ -\!1\!\leq\!w\!<\!v\!<\!1$&$c_{15}^1\!=\!c_{25}^2\!=\!1,\ c_{35}^3\!=\!v,\ c_{45}^4\!=\!w$&&$p_1$&&$(S_{12},a,b,1)$\\[0.9ex]
$A_{5,7}^{u,u.u}$&$u\!\neq\!0,\ -\!1\!\leq\!u\!<\!1$&$c_{15}^1\!=\!1,\ c_{25}^2\!=\!c_{35}^3\!=\!c_{45}^4\!=\!u$&&$p_1,p_2$&$E_2^2$&$(1,S_{234},1)$\\[0.9ex]
$A_{5,7}^{1,1,w}$&$w\!\neq\!0,\ -\!1\!\leq\!w\!<\!1$&$c_{15}^1\!=\!c_{25}^2\!=\!c_{35}^3\!=\!1,\ c_{45}^4\!=\!w$&&$p_1$&&$(S_{123},a,1)$\\[0.9ex]
$A_{5,7}^{1,1,1}$&&$c_{15}^1\!=\!c_{25}^2\!=\!c_{35}^3\!=\!c_{45}^4\!=\!1$&&$p_1$&&$(S_{1234},1)$\\[0.9ex]
$A_{5,7}^{1,v,v}$&$0\!<\!|v|\!<\!1$&$c_{15}^1\!=\!c_{25}^2\!=\!1,\ c_{35}^3\!=\!c_{45}^4\!=\!v$&&$p_1,p_3$&$E_3^3$&$(S_{12},S_{34},1)$\\[0.9ex]
$A_{5,7}^{u,v,-\!1}$&$uv(u\!+\!v)\!\neq\!0,\ -\!1\!<\!v\!<\!u\!<\!1$&$c_{15}^1\!=\!1,\ c_{25}^2\!=\!u,\ c_{35}^3\!=\!v,\ c_{45}^4\!=\!-\!1$&&$p_1$&&$(1,a,b,c,1)$\\[0.9ex]
$A_{5,7}^{u,-\!u,-\!1}$&$0\!<\!u\!<\!1$&$c_{15}^1\!=\!1,\ c_{25}^2\!=\!u,\ c_{35}^3\!=\!-u,\ c_{45}^4\!=\!-\!1$&&$(-\!X_4,X_3,X_2,X_1,-\!X_5)$&&$(1,a,b,c,1)$\\[0.9ex]
$A_{5,7}^{1,-\!1,-\!1}$&&$c_{15}^1\!=\!c_{25}^2\!=\!1,\ c_{35}^3\!=\!c_{45}^4\!=\!-\!1$&&$p_3,\ (-\!X_4,X_3,X_2,X_1,-\!X_5)$&$E_3^3$&$(S_{12},S_{34},1)$\\[0.9ex]
$A_{5,8}^{u}$&$0\!<\!|u|\!<\!1$&$c_{25}^1\!=\!c_{35}^3\!=\!1,\ c_{45}^4\!=\!u$&&&$E_5^2$&$(a,a,b,c,1)$\\[0.9ex]
$A_{5,8}^{-\!1}$&&$c_{25}^1\!=\!c_{35}^3\!=\!1,\ c_{45}^4\!=\!-\!1$&&$(X_1,-\!X_2,-\!X_4,X_3,-\!X_5)$&$E_3^3\!+\!E_4^4,\ E_5^2$&$(a,a,b,1,1)$\\[0.9ex]
$A_{5,8}^{1}$&&$c_{25}^1\!=\!c_{35}^3\!=\!c_{45}^4\!=\!1$&&$p_3$&$E_3^3,\ E_5^2$&$(a,a,S_{34},1)$\\[0.9ex]
$A_{5,9}^{u,v}$&$v\!<\!u,\ u,v\not\in\{0,1\}$&$c_{15}^1\!=\!c_{25}^1\!=\!c_{25}^2\!=\!1,\ c_{35}^3\!=\!u,\
c_{45}^4\!=\!v$&&&&$(a,a,b,c,1)$\\[0.9ex]
$A_{5,9}^{1,v}$&$v\!\neq\!0,\ v\!<\!1$&$c_{15}^1\!=\!c_{25}^1\!=\!c_{25}^2\!=\!c_{35}^3\!=\!1,\
c_{45}^4\!=\!v$&&&$E_2^3,\ E_3^1$&$(a,a,b,c,1)$\\[0.9ex]
$A_{5,9}^{u,1}$&$1\!<\!u$&$c_{15}^1\!=\!c_{25}^1\!=\!c_{25}^2\!=\!c_{45}^4\!=\!1,\ c_{35}^3\!=\!u$&&&$E_2^4,\ E_4^1$&$(a,a,b,c,1)$\\[0.9ex]
$A_{5,9}^{u,u}$&$u\!\not\in\!\{0,1\}$&$c_{15}^1\!=\!c_{25}^1\!=\!c_{25}^2\!=\!1,\
c_{35}^3\!=\!c_{45}^4\!=\!u$&&$p_3$&$E_3^3$&$(a,b,S_{34},1)$\\[0.9ex]
$A_{5,9}^{1,1}$&&$c_{15}^1\!=\!c_{25}^1\!=\!c_{25}^2\!=\!c_{35}^3\!=\!c_{45}^4\!=\!1$&&$p_3$&$E_2^3,\
E_2^4,\ E_3^1,\ E_3^3,\ E_4^1$&$(a,b,S_{34},1)$\\[0.9ex]
$A_{5,10}$&&$c_{25}^1\!=\!c_{35}^2\!=\!c_{45}^4\!=\!1$&&&$E_3^1,\ E_5^3$&$(a,a,a,b,1)$\\[0.9ex]
$A_{5,11}^{u}$&$u\!\not\in\!\{0,1\}$&$c_{15}^1\!=\!c_{25}^1\!=\!c_{25}^2\!=\!c_{35}^2\!=\!c_{35}^3\!=\!1,\
c_{45}^4\!=\!u$&&&$E_3^1$&$(a,a,a,b,1)$\\[0.9ex]
$A_{5,11}^{1}$&&$c_{15}^1\!=\!c_{25}^1\!=\!c_{25}^2\!=\!c_{35}^2\!=\!c_{35}^3\!=\!c_{45}^4\!=\!1$&&&$E_3^1,\
E_3^4,\ E_4^1$&$(a,a,a,b,1)$\\[0.9ex]
$A_{5,12}$&&$c_{15}^1\!=\!c_{25}^1\!=\!c_{25}^2\!=\!c_{35}^2\!=\!c_{35}^3\!=\!c_{45}^3\!=\!c_{45}^4\!=\!1$&&&$E_3^1\!+\!E_4^2,\
E_4^1$&$(a,a,a,a,1)$\\[0.9ex]
$A_{5,13}^{u,v,w}$&$uvw\!\neq\!0,\ -\!1\!\leq\!u\!<\!1,$ & \parbox{4cm}{\rule{0pt}{\baselineskip}$\!c_{15}^1\!=\!1,\ c_{25}^2\!=\!u,\ c_{35}^3\!=\!c_{45}^4\!=\!v,$\\[2pt]$c_{35}^4\!=\!-\!w,\ c_{45}^3\!=\!w$\vspace{1ex}} &&&$\big[ E_3^3\!+\!E_4^4 \big]_{v\!/\!w}$&$\!(a,b,1,1,1)$\\[0.9ex]
$A_{5,13}^{u,0,w}$&$uw\!\neq\!0,\ \left|u\right|\!<\!1,$&$\!c_{15}^1\!=\!1,\ c_{25}^2\!=\!u,\ c_{35}^4\!=\!-\!w,\ c_{45}^3\!=\!w$&&&$E_3^3\!+\!E_4^4$&$\!(a,b,1,1,1)$\\[0.9ex]
$A_{5,13}^{1,v,w}$&$vw\neq 0$ & \parbox{4cm}{\rule{0pt}{1\baselineskip}$c_{15}^1\!=\!c_{25}^2\!=\!1,\ c_{35}^3\!=\!c_{45}^4\!=\!v,$\\[2pt]$c_{35}^4\!=\!-\!w,\ c_{45}^3\!=\!w$\vspace{1ex}}&&$p_1$&$E_1^1,\ \big[ E_3^3\!+\!E_4^4 \big]_{v\!/\!w}$&$(S_{12},1,1,1)$\\[0.9ex]
$A_{5,13}^{1,0,w}$&$w\neq 0$&$c_{15}^1\!=\!c_{25}^2\!=\!1,\ c_{35}^4\!=\!-\!w,\ c_{45}^3\!=\!w$&&$p_1$&$E_1^1,\ E_3^3\!+\!E_4^4$&$(S_{12},1,1,1)$\\[0.9ex]
$A_{5,13}^{-\!1,0,w}$&$w\!\neq\!0$&$c_{15}^1\!=\!1,\ c_{25}^2\!=\!-\!1,\ c_{35}^4\!=\!-\!w,\
c_{45}^3\!=\!w$&&$(-\!X_2,X_1,X_3,-\!X_4,-\!X_5)$&$E_1^1\!+\!E_2^2,\ E_3^3\!+\!E_4^4$&$\!(a,1,1,1,1)$\\[0.9ex]
$A_{5,14}^u$&$u\neq 0$&$c_{25}^1\!=\!c_{45}^3\!=\!1,\ c_{35}^3\!=\!c_{45}^4\!=\!u,\
c_{35}^4\!=-\!1$&&&$E_2^1,\ \big[ E_3^3\!+\!E_4^4 \big]_u,\ E_5^2$&$(a,a,1,1,1)$\\[0.9ex]
$A_{5,14}^0$&&$c_{25}^1\!=\!c_{45}^3\!=\!1, \ c_{35}^4\!=\!-\!1$&&$p_{245}$&$E_2^1,\ E_3^3\!+\!E_4^4,\ E_5^2$&$(a,a,1,1,1)$\\[0.9ex]
$A_{5,15}^u$&$0< \left|u\right|<1$&$c_{15}^1\!=\!c_{25}^1\!=\!c_{25}^2\!=\!c_{45}^3\!=\!1,\ c_{35}^3\!=\!c_{45}^4\!=\!u$&&&$E_2^1$&$(a,a,b,b,1)$\\[0.9ex]
$A_{5,15}^{-\!1}$&&$c_{15}^1\!=\!c_{25}^1\!=\!c_{25}^2\!=\!c_{45}^3\!=\!1,\
c_{35}^3\!=\!c_{45}^4\!=\!-\!1$&&$(-\!X_3,X_4,X_1,-\!X_2,-\!X_5)$&$E_2^1,\ E_3^3\!+\!E_4^4$&$(a,a,1,1,1)$\\[0.9ex]
$A_{5,15}^0$&&$c_{15}^1\!=\!c_{25}^1\!=\!c_{25}^2\!=\!c_{45}^3\!=\!1$&&&$E_2^1,\ E_5^4$&$(a,a,b,b,1)$\\[0.9ex]
$A_{5,15}^1$(*)&&$c_{15}^1\!=\!c_{25}^2\!=\!c_{35}^1\!=\!c_{35}^3\!=\!c_{45}^2\!=\!c_{45}^4=1$&&$p_{13}$&$E_3^1,\ E_3^2,\ E_4^1,\ E_4^2$&$(S_{12},S_{12},1)$\\[0.9ex]
$A_{5,16}^{u,v}$&$v\ne 0$ & \parbox{5.5cm}{\rule{0pt}{1\baselineskip}$c_{15}^1\!=\!c_{25}^1\!=\!c_{25}^2\!=1,\ c_{35}^3\!=\!c_{45}^4\!=u,$\\[2pt]$\!c_{35}^4\!=\!-\!v,\ c_{45}^3\!=v$\vspace{1ex}} &&& $E_2^1,\ \big[ E_3^3\!+\!E_4^4 \big]_{u\!/\!v}$ &$(a,a,1,1,1)$\\[0.9ex]
$A_{5,17}^{w,u,v}$&$w\neq 0$ &  \parbox{5.5cm}{\rule{0pt}{1\baselineskip}$c_{15}^1\!=\!c_{25}^2\!=u,\ c_{35}^3\!=\!c_{45}^4\!=v,$\\[2pt]$c_{15}^2\!=\!-\!1,\ c_{25}^1\!=\!1,\ c_{35}^4\!=\!-w,\ c_{45}^3\!=w$\vspace{1ex}}&&&&see \S4.2\\[0.9ex]
$A_{5,18}^u$&$u>0$ & \parbox{5.5cm}{\rule{0pt}{1\baselineskip}$c_{15}^1\!=\!c_{25}^2\!=\!c_{35}^3\!=\!c_{45}^4\!=u,$\\[2pt]$c_{15}^2\!=\!c_{35}^4\!=\!-\!1,\ c_{25}^1\!=\!c_{35}^1\!=\!c_{45}^2\!=\!c_{45}^3\!=\!1$\vspace{1ex}}&&&
\begin{minipage}[c]{3.8cm}\vspace{1ex}\begin{center}$\big[ E_1^1\!+\!E_2^2\!+\!E_3^3\!+\!E_4^4 \big]_u$,\\[2pt]$E_3^1\!+\!E_4^2,\
E_3^2\!-\!E_4^1$\vspace{1ex}\end{center}\end{minipage}
&$\!(1,1,1,1,1)$\\[0.9ex]
$A_{5,18}^0$&&$c_{15}^2\!=\!c_{35}^4\!=\!-\!1,c_{25}^1\!=\!c_{35}^1\!=\!c_{45}^2\!=\!c_{45}^3\!=\!1$&&$p_{235}$&\begin{minipage}[c]{3.8cm}\begin{center}\vspace{1ex}$E_1^1\!+\!E_2^2\!+\!E_3^3\!+\!E_4^4$,\\[2pt]$E_3^1\!+\!E_4^2,\ E_3^2\!-\!E_4^1$\vspace{1ex}\end{center}\end{minipage}&$\!(1,1,1,1,1)$\\[0.9ex]
$A_{5,19}^{u,v}$ & $u\not\in\{0,1,2\},\ v\not\in\{0,u\!-\!1,u\}$ & $c_{15}^1\!=\!u,\ c_{23}^1\!=\!c_{25}^2\!=\!1,\
c_{35}^3\!=\!u\!-\!1,\ c_{45}^4\!=\!v$&&$p_{12}$&&$(a,1,a,b,1)$\\[0.9ex]
$A_{5,19}^{u,u}$&$u\notin \{0,1,2\}$&$c_{23}^1\!=\!c_{25}^2\!=\!1,\ c_{35}^3\!=\!u\!-\!1,\ c_{15}^1\!=\!c_{45}^4\!=\!u$&&$p_{12}$&$E_4^1$&$(a,1,a,b,1)$\\[0.9ex]
$A_{5,19}^{u,u\!-\!1}$&$u\notin \{0,1,2\}$&$c_{15}^1\!=\!u,\ c_{23}^1\!=\!c_{25}^2\!=\!1,\ c_{35}^3\!=\!c_{45}^4\!=\!u\!-\!1$&&$p_{12}$&$E_3^4$&$(a,1,a,b,1)$\\[0.9ex]
$A_{5,19}^{0,v}$&$v\notin\{-1,0,1\}$&$c_{23}^1\!=\!c_{25}^2\!=\!1,\ c_{35}^3\!=\!-\!1,\ c_{45}^4\!=\!v$&&$p_{12}$&$E_5^1$&$(a,1,a,b,1)$\\[0.9ex]
$A_{5,19}^{0,-\!1}$&&$c_{23}^1\!=\!c_{25}^2\!=\!1,\ c_{35}^3\!=\!c_{45}^4\!=\!-\!1$&&$p_{12}$&$E_3^4,\
E_5^1$&$(a,1,a,b,1)$\\[0.9ex]
$A_{5,19}^{0,1}$&&$c_{23}^1\!=\!c_{25}^2\!=\!c_{45}^4\!=\!1,\ c_{35}^3\!=\!-\!1$&&$p_{12}$&$E_2^4,\ E_5^1$&$(a,1,a,b,1)$\\[0.9ex]
$A_{5,19}^{1,v}$&$v\notin\{0,1\}$&$c_{15}^1\!=\!c_{23}^1\!=\!c_{25}^2\!=\!1,\ c_{45}^4\!=\!v$&&$p_{12}$&&$(a,1,a,b,1)$\\[0.9ex]
$A_{5,19}^{1,1}$&&$c_{15}^1\!=\!c_{23}^1\!=\!c_{25}^2\!=\!=\!c_{45}^4\!1$&&$p_{12}$&$E_2^4,\ E_4^1$&$(a,1,a,b,1)$\\[0.9ex]
$A_{5,19}^{2,v}$&$v\notin\{0,1,2\}$&$c_{15}^1\!=\!2,\ c_{23}^1\!=\!c_{25}^2\!=\!c_{35}^3\!=\!1,\ c_{45}^4\!=\!v$&&$p_{12}$&&$(1,S_{23},a,1)$\\[0.9ex]
$A_{5,19}^{2,1}$&&$c_{23}^1\!=\!c_{25}^2=c_{35}^3\!=\!c_{45}^4\!=\!1,\
c_{15}^1\!=\!2$&&$p_{12}$&$E_2^4,\ E_3^4$&$(1,S_{23},a,1)$\\[0.9ex]
$A_{5,19}^{2,2}$&&$c_{23}^1\!=\!c_{25}^2\!=\!c_{35}^3\!=\!1,\ c_{15}^1\!=\!c_{45}^4\!=\!2$&&$p_{12}$&$E_4^1$&$(1,S_{23},a,1)$\\[0.9ex]
$A_{5,20}^u$&$u\notin \{0,1,2\}$ &  \parbox{5.5cm}{\rule{0pt}{1\baselineskip}$c_{15}^1\!=\!c_{45}^4\!=\!u,\ c_{23}^1\!=\!c_{25}^2\!=\!c_{45}^1\!=\!1,$\\[2pt]$c_{35}^3\!=\!u\!-\!1$\vspace{1ex}}&&&&$(ab,a,b,ab,1)$\\[0.9ex]
$A_{5,20}^0$&&$c_{23}^1\!=\!c_{25}^2\!=\!c_{45}^1\!=\!1,\ c_{35}^3\!=\!-\!1$&&$(X_1,X_3,-\!X_2,-\!X_4,-\!X_5)$&$E_2^2\!-\!E_3^3,\ E_5^4$&$(a,1,a,a,1)$\\[0.9ex]
$A_{5,20}^1$&&$c_{15}^1\!=\!c_{23}^1\!=\!c_{25}^2\!=\!c_{45}^1\!=\!c_{45}^4\!=\!1$&&&$E_5^3\!-\!E_2^4$&$(ab,a,b,ab,1)$\\[0.9ex]
$A_{5,20}^2$&&$c_{15}^1\!=\!c_{45}^4\!=\!2,\
c_{23}^1\!=\!c_{25}^2\!=\!c_{35}^3\!=\!c_{45}^1\!=\!1$&&$p_{124}$&$E_1^1\!+\!E_2^2\!+\!E_4^4$&$(1,S_{23},1,1)$\\[0.9ex]
$A_{5,21}$&&$c_{15}^1\!=\!2,\
c_{23}^1\!=\!c_{25}^2\!=\!c_{25}^3\!=\!c_{35}^3\!=\!c_{35}^4\!=\!c_{45}^4\!=\!1$&&&$E_{2}^4$&$(a^2\!,a,a,a,1)$\\[0.9ex]
$A_{5,22}$&&$c_{23}^1\!=\!c_{25}^3\!=\!c_{45}^4\!=\!1$&&&$E_5^1$&$(a^2\!,a,a,b,1)$\\[0.9ex]
$A_{5,23}^u$&$u\notin \{0,1,2\}$&$c_{15}^1\!=\!2,\
c_{23}^1\!=\!c_{25}^2\!=\!c_{25}^3\!=\!c_{35}^3\!=\!1,\ c_{45}^4\!=\!u$&&&&$(a^2\!,a,a,b,1)$\\[0.9ex]
$A_{5,23}^1$&&$c_{15}^1\!=\!2,\
c_{23}^1\!=\!c_{25}^2\!=\!c_{25}^3\!=\!c_{35}^3\!=\!c_{45}^4\!=\!1$&&&$E_2^4$&$(a^2\!,a,a,b,1)$\\[0.9ex]
$A_{5,23}^2$&&$c_{15}^1\!=\!c_{45}^4\!=\!2,\
c_{23}^1\!=\!c_{25}^2\!=\!c_{25}^3\!=\!c_{35}^3\!=\!1$&&&$E_4^1$&$(a^2\!,a,a,b,1)$\\[0.9ex]
$A_{5,24}^\epsilon$&$\epsilon=\pm 1$ &  \parbox{5.5cm}{\rule{0pt}{1\baselineskip}$c_{15}^1\!=\!c_{45}^4\!=\!2,\ c_{23}^1\!=\!c_{25}^2\!=\!c_{25}^3\!=\!c_{35}^3\!=\!1,$\\[2pt]$c_{45}^1\!=\!\epsilon$\vspace{1ex}}&&&$E_4^1$&$(a^2\!,a,a,a^2\!,1)$\\[0.9ex]
$A_{5,25}^{u,v}$&$u\notin \{0,2v\},\ v\neq 0$& \parbox{5.5cm}{\rule{0pt}{1\baselineskip}$c_{15}^1\!=\!2v,\ c_{23}^1\!=\!c_{25}^3\!=\!1,\ c_{25}^2\!=\!c_{35}^3\!=\!v,$\\[2pt]$c_{35}^2\!=\!-\!1,\ c_{45}^4\!=\!u$\vspace{1ex}}&&&$\big[ 2E_1^1\!+\!E_2^2\!+\!E_3^3 \big]_v$&$(1,1,1,a,1)$\\[0.9ex]
$A_{5,25}^{u,0}$&$u\neq 0$&$c_{23}^1\!=\!c_{25}^3\!=\!1,\ c_{35}^2\!=\!-\!1,\
c_{45}^4\!=\!u$&&&$2E_1^1\!+\!E_2^2\!+\!E_3^3,\ E_5^1$&$(1,1,1,a,1)$\\[0.9ex]
$A_{5,25}^{2v,v}$&$v\neq 0$& \parbox{5.5cm}{\rule{0pt}{1\baselineskip}$c_{15}^1\!=\!c_{45}^4\!=\!2v,\ c_{23}^1\!=\!c_{25}^3\!=\!1,\
c_{25}^2\!=\!c_{35}^3\!=\!v,$\\[2pt]$c_{35}^2\!=\!-1$\vspace{1ex}}&&&$\big[ 2E_1^1\!+\!E_2^2\!+\!E_3^3 \big]_v,\ E_4^1$&$(1,1,1,a,1)$\\[0.9ex]
$A_{5,26}^{u,\epsilon}$&$\epsilon=\pm1,\ u\neq 0$& \parbox{5.5cm}{\rule{0pt}{1\baselineskip}$c_{15}^1\!=\!c_{45}^4\!=\!2u,\ c_{23}^1\!=\!c_{25}^3\!=\!1,\ c_{25}^2\!=\!c_{35}^3\!=\!u,$\\[2pt]$c_{35}^2\!=\!-1,c_{45}^1\!=\!\epsilon$\vspace{1ex}}&&&\begin{minipage}[c]{3.8cm}\begin{center}\vspace{1ex}$\big[ 2E_1^1\!+\!E_2^2\!+\!E_3^3\!+\!2E_4^4 \big]_u$,\\[2pt]$E_4^1$\vspace{1ex}\end{center}\end{minipage}&$(1,1,1,1,1)$\\[0.9ex]
$A_{5,26}^{0,\epsilon}$&$\epsilon=\pm 1$&$c_{23}^1\!=\!c_{25}^3\!=\!1,\ c_{35}^2\!=\!-\!1,\
c_{45}^1\!=\!\epsilon$&&$p_{135}$&\begin{minipage}[c]{3.8cm}\begin{center}\vspace{1ex}$2E_1^1\!+\!E_2^2\!+\!E_3^3\!+\!2E_4^4$,\\[2pt]$E_4^1,\
E_5^4$\vspace{1ex}\end{center}\end{minipage}&$(1,1,1,1,1)$\\[0.9ex]
$A_{5,27}$&&$c_{15}^1\!=\!c_{23}^1\!=\!c_{35}^3\!=\!c_{35}^4\!=\!c_{45}^1\!=\!c_{45}^4\!=\!1$&&&$E_5^2\!-\!E_4^1$&$(a,1,a,a,1)$\\[0.9ex]
$A_{5,28}^u$&$u\notin\{0,1,2\}$& \parbox{5.5cm}{\rule{0pt}{1\baselineskip}$c_{15}^1\!=\!u,\ c_{23}^1\!=\!c_{35}^3\!=\!c_{35}^4\!=\!c_{45}^4\!=\!1,$\\[2pt]$c_{25}^2\!=\!u\!-\!1$\vspace{1ex}}&&&&$(ab,a,b,b,1)$\\[0.9ex]
$A_{5,28}^0$&&$c_{23}^1\!=\!c_{35}^3\!=\!c_{35}^4\!=\!c_{45}^4\!=\!1,\ c_{25}^2\!=\!-\!1$&&&$E_5^1$&$(ab,a,b,b,1)$\\[0.9ex]
$A_{5,28}^1$&&$c_{15}^1\!=\!c_{23}^1\!=\!c_{35}^3\!=\!c_{35}^4\!=\!c_{45}^4\!=\!1$&&&$E_5^2\!-\!E_4^1$&$(ab,a,b,b,1)$\\[0.9ex]
$A_{5,28}^2$&&$c_{15}^1\!=\!2,\
c_{23}^1\!=\!c_{25}^2\!=\!c_{35}^3\!=\!c_{35}^4\!=\!c_{45}^4\!=\!1$&&&$E_2^4,\ E_3^2$&$(ab,a,b,b,1)$\\[0.9ex]
$A_{5,29}$&&$c_{15}^1\!=\!c_{24}^1\!=\!c_{25}^2\!=\!c_{45}^3\!=\!1$&&&$E_{2}^1$&$(ab,a,b,b,1)$\\[0.9ex]
$A_{5,30}^u$&$u\notin\{-1,2\}$& \parbox{5.5cm}{\rule{0pt}{1\baselineskip}$c_{15}^1\!=\!u\!+\!1,\ c_{24}^1\!=\!c_{34}^2\!=\!c_{45}^4\!=\!1,\ c_{25}^2\!=\!u,$\\[2pt]$c_{35}^3\!=\!u\!-\!1$\vspace{1ex}}&&$p_{24}$&&$(a,a,a,1,1)$\\[0.9ex]
$A_{5,30}^{-\!1}$&&$c_{24}^1\!=\!c_{34}^2\!=\!c_{45}^4\!=\!1,\ c_{25}^2\!=\!-\!1,\ c_{35}^3\!=\!-\!2$&&$p_{24}$&$E_5^1$&$(a,a,a,1,1)$\\[0.9ex]
$A_{5,30}^2$&&$c_{24}^1\!=\!c_{34}^2\!=\!c_{45}^4\!=\!c_{35}^3\!=\!1,\ c_{25}^2\!=\!2,\ c_{15}^1\!=\!3$&&$p_{24}$&$E_4^3$&$(a,a,a,1,1)$\\[0.9ex]
$A_{5,31}$&& \parbox{5.5cm}{\rule{0pt}{1\baselineskip}$c_{15}^1\!=\!3,\ c_{24}^1\!=\!c_{34}^2\!=\!c_{35}^3\!=\!c_{45}^3\!=\!c_{45}^4\!=\!1,$\\[2pt]$c_{25}^2\!=\!2$\vspace{1ex}}&&&&$(a^3\!,a^2\!,a,a,1)$\\[0.9ex]
$A_{5,32}^u$&$u\neq 0$&$c_{15}^1\!=\!c_{24}^1\!=\!c_{25}^2\!=\!c_{34}^2\!=\!c_{35}^3\!=\!1,\ c_{35}^1\!=\!u$&&$p_{24}$&&$(a,a,a,1,1)$\\[0.9ex]
$A_{5,32}^0$&&$c_{15}^1\!=\!c_{24}^1\!=\!c_{25}^2\!=\!c_{34}^2\!=\!c_{35}^3\!=\!1$&&$p_{123}$&$E_3^1$&$(a^2\!,a,1,a.1)$\\[0.9ex]
$A_{5,33}^{u,v}$&\parbox{3.5cm}{\rule{0pt}{1\baselineskip}$u^2\!+\!v^2\!\neq\! 0,\ u\!\neq\!-\!1\!\neq v,$\\[2pt]$v\!<\!u,\ (u,\!v)\!\neq\!(1,\!0)$\vspace{1ex}} &$c_{14}^1\!=\!c_{25}^2\!=\!1,\ c_{34}^3\!=\!v,\ c_{35}^3\!=\!u$&&$p_1,p_2$&&$(1,1,a,1,1)$\\[0.9ex]
$A_{5,33}^{u,u}$&$u\not\in\{\!-\!1,\!0\}$&$c_{14}^1\!=\!c_{25}^2\!=\!1,\ c_{34}^3\!=\!c_{35}^3\!=\!u$&&$p_2,\
(X_2,-\!X_1,X_3,X_5,X_4)$&&$(1,1,a,1,1)$\\[0.9ex]
$A_{5,33}^{u,-\!1}$&$u\neq -\!1$&$c_{14}^1\!=\!c_{25}^2\!=\!1,\ c_{34}^3\!=\!-\!1,\ c_{35}^3\!=\!u$&&$p_2,\ (-\!X_3,X_2,X_1,\!-\!X_4,X_5\!+\!uX_4)$&&$(1,1,a,1,1)$\\[0.9ex]
$A_{5,33}^{-\!1,-\!1}$&& $c_{14}^1\!=\!c_{25}^2\!=\!1,c_{34}^3\!=\!c_{35}^3\!=\!-\!1$ & & \begin{minipage}[c]{3.8cm}\begin{center}\vspace{1ex}$(X_3,-\!X_1,X_2,-\!X_5,X_4\!-\!X_5),$\\[2pt]$(X_2,-\!X_1,X_3,X_5,X_4)$\vspace{1ex}\end{center}\end{minipage}&$E_3^3$&$(1,1,1,1,1)$\\[0.9ex]
$A_{5,33}^{1,0}$&&$c_{14}^1\!=\!c_{25}^2\!=\!c_{35}^3\!=\!1$&&$p_1,p_2$&&$(1,S_{23},1,1)$\\[0.9ex]
$A_{5,34}^u$&&$c_{14}^1\!=\!u,\ c_{15}^1\!=\!c_{24}^2\!=\!c_{34}^3\!=\!c_{35}^2\!=\!1$&&$p_{23}$&&$(a,1,1,1,1)$\\[0.9ex]
$A_{5,35}^{u,v}$&$u\ne 0$& \parbox{5.5cm}{\rule{0pt}{1\baselineskip}$c_{14}^1\!=\!v,\ c_{15}^1\!=\!u,\ c_{24}^2\!=\!c_{34}^3\!=\!c_{35}^2\!=\!1,$\\[2pt]$c_{25}^3\!=\!-1$\vspace{1ex}}&&&&$(a,1,1,1,1)$\\[0.9ex]
$A_{5,35}^{0,v}$&&$c_{14}^1\!=\!v,\ c_{24}^2\!=\!c_{34}^3\!=\!c_{35}^2\!=\!1,\ c_{25}^3\!=\!-\!1$&&$p_{35}$&&$(a,1,1,1,1)$\\[0.9ex]
$A_{5,36}$&&$c_{14}^1\!=\!c_{23}^1\!=\!c_{24}^2\!=\!c_{35}^3\!=\!1,\ c_{25}^2\!=\!-\!1$&&$p_{13},\ (X_1,X_3,\!-\!X_2,X_4\!+\!X_5,\!-\!X_5)$&&$(1,1,1,1,1)$\\[0.9ex]
$A_{5,37}$&& \parbox{5.5cm}{\rule{0pt}{1\baselineskip}$c_{14}^1\!=\!2,\ c_{23}^1\!=\!c_{24}^2\!=\!c_{34}^3\!=\!c_{35}^2\!=\!1,$\\[2pt]$c_{25}^3\!=\!-1$\vspace{1ex}}&&$p_{135}$&&$(1,1,1,1,1)$\\[0.9ex]
$A_{5,38}$&&$c_{14}^1\!=\!c_{25}^2\!=\!c_{45}^3\!=\!1$&&$(-X_2,X_1,-\!X_3,X_5,X_4)$&$E_2^2$&$(a,1,1,1,1)$\\[0.9ex]
$A_{5,39}$&&$c_{14}^1\!=\!c_{24}^2\!=\!c_{25}^1\!=\!c_{45}^3\!=\!1,\ c_{15}^2\!=\!-\!1$&&$p_{135}$&$E_1^1\!+\!E_2^2,\ E_4^3$&$(1,1,1,1,1)$\\[0.9ex]
$A_{5,40}$&$\mathfrak{sl}(2,\mathbb{R})\ltimes\mathbb{R}^2$& \parbox{5.5cm}{\rule{0pt}{1\baselineskip}$c_{12}^1\!=\!c_{23}^3\!=\!2,\ c_{13}^2\!=\!c_{25}^5\!=\!-\!1,$\\[2pt]$c_{14}^5\!=\!c_{24}^4\!=\!c_{35}^4\!=\!1$\vspace{1ex}}&&$p_{135},\ ((-\!X_3,-\!X_2,-\!X_1,-\!X_5,X_4))$&$E_4^4\!+\!E_5^5$&$(1,1,1,1,1)$\\[0.9ex]
\end{longtable}
\end{center}
\end{landscape}

\restoregeometry 

\normalsize

\section*{Acknowledgements}

We thank Alessandro Torrielli for his helpful comments on an early draft, and Katy Pelling, whose project work on low-dimensional real Lie algebras contributed to a concise presentation of the automorphisms. We also thank the Nuffield Foundation for an Undergraduate Research Bursary (URB/01451/G).

\end{document}